\documentstyle[12pt]{article}

\def\R{\hat{R}}

\def\lb{\label}
\def\be{\begin{equation}}
\def\ee{\end{equation}}

\begin{document}

\begin{center}
\Large{ STANDARD COMPLEX FOR QUANTUM\\
LIE ALGEBRAS}
\end{center}

\vspace{.5cm}

\begin{center}
\large{C.\,BURDIK${}^{a}$, A.P.\,ISAEV${}^{b}$ 
and O.\,OGIEVETSKY${}^{c}$ }
\end{center}

\begin{center}
${}^{a}$ Department of Mathematics, Faculty of Nuclear Sciences
and Physical Engineering, Trojanova 13, 120 00 Prague 2,
Czech Republic \\
E-mail: burdik@dec1.fjfi.cvut.cz
\\ 
\vspace{0.3cm}
${}^{b}$ Bogoliubov Theoretical Laboratory,
Joint Institute for Nuclear Research,
Dubna, Moscow region 141980, Russia \\
E-mail: isaevap@thsun1.jinr.ru
\\ 
\vspace{0.3cm}
${}^{c}$ Center of Theoretical Physics, Luminy,
13288 Marseille, France \\ 
and P. N. Lebedev Physical Institute, Theoretical Department,
Leninsky pr. 53, 117924 Moscow, Russia \\
E-mail: oleg@cpt.univ-mrs.fr
\end{center}

\vspace{1cm}

\begin{center}
ABSTRACT
\end{center}
\noindent
For a quantum Lie algebra $\Gamma$, let $\Gamma^\wedge$  be
its exterior extension (the algebra $\Gamma^\wedge$
is canonically defined).
We introduce a differential on the exterior extension algebra $\Gamma^\wedge$
which provides the structure of a complex on $\Gamma^{\wedge}$. 
In the situation when $\Gamma$ is a usual Lie algebra this
complex coincides with the "standard complex".
The differential  is realized as a commutator with
a (BRST) operator $Q$ in a larger algebra $\Gamma^\wedge[\Omega]$, with extra 
generators canonically conjugated to the exterior generators
of $\Gamma^{\wedge}$. A recurrent relation which defines uniquely
the operator $Q$ is given.

\pagebreak 

{\bf 1.} A quantum Lie
algebra \cite{Wor1}, \cite{Wor},
\cite{Ber}, \cite{Cast} is defined by two tensors
$C^k_{ij}$ and $\sigma^{mk}_{ij}$ (indices belong
to some set ${\cal N}$, say, ${\cal N}= \{ 1, \dots, N \}$). 
By definition, the matrix $\sigma^{mk}_{ij}$ has an eigenvalue
1; one demands that $(P_{(1)})^{mk}_{ij} \, C^n_{mk} = 0$, where
$P_{(1)}$ is a projector on the eigenspace of $\sigma$ corresponding to the
eigenvalue 1. 

By definition, a quantum Lie algebra $\Gamma$
is generated by elements $\chi_i$, $i=1, \dots, N$, subjected
to relations
 \be
 \lb{int1}
\chi_i \, \chi_j - \sigma^{mk}_{ij} \chi_m \, \chi_k = C^k_{ij} \,
\chi_k \; .
 \ee
Here the structure constants $C^k_{ij}$ obey
 \be
 \lb{int2}
 \begin{array}{c}
C^p_{ni} \, C^{l}_{pj} = \sigma^{mk}_{ij} \, C^p_{nm} \,
C^{l}_{pk} + C^p_{ij} \, C^{l}_{np} \; \Leftrightarrow \\ \\
\Leftrightarrow \;
C^{<1|}_{|1 2>} \, C^{<4|}_{|1 3>} = \sigma_{23} \, C^{<1|}_{|1 2>} \,
C^{<4|}_{|1 3>} + C^{<3|}_{|23>} \, C^{<4|}_{|1 3>} \; ,
\end{array}
 \ee
 \be
 \lb{int3}
C^k_{ni} \, \sigma^{pm}_{kq} =  \sigma^{sj}_{iq} \,  
\sigma^{pk}_{ns} \, C^m_{kj} 
\;\;\; \Leftrightarrow \;\;\;
C^{<1|}_{|1 2>} \, \sigma_{13} = 
\sigma_{23} \, \sigma_{12} \, C^{<3|}_{|2 3>} \; ,
 \ee
 
 \be
 \lb{int3a}
 \begin{array}{c}
(\sigma^{pj}_{im} \, C^n_{qp} + \delta^n_q \, C^j_{im} ) \,
\sigma^{ks}_{nj}
= \sigma^{jn}_{qi} \, (\sigma^{ps}_{nm} \, C^k_{jp} +
\delta^k_j  \, C^{s}_{nm})  \;\;\; \Leftrightarrow  \\ \\
 \Leftrightarrow \;\;\;
(\sigma_{23} \, C^{<1|}_{|12>} +  C^{<3|}_{|23>} ) \,
\sigma_{13} = \sigma_{12} \, (\sigma_{23} \, C^{<1|}_{|12>} +
C^{<3|}_{|23>} ) \; .
\end{array}
 \ee
The matrix $\sigma^{mk}_{ij}$ satisfies the Yang-Baxter equation
 \be
 \lb{int1a}
\sigma^{j_1 j_2}_{i_1 i_2} \, \sigma^{n_2 k_3}_{j_2 i_3} \,
\sigma^{k_1 k_2}_{j_1 n_2} = \sigma^{j_2 j_3}_{i_2 i_3} \,
\sigma^{k_1 n_2}_{i_1 j_2} \, \sigma^{k_2 k_3}_{n_2 j_3} \;\;\;
\Leftrightarrow \;\;\;  \sigma_{12} \, \sigma_{23} \, \sigma_{12} =
\sigma_{23} \, \sigma_{12} \, \sigma_{23} \; .
 \ee
In the right hand 
side of (\ref{int3})-(\ref{int1a}) we
use FRT matrix notations \cite{FRT}; $\{ 1,2,3, \dots \}$
are the numbers of vector spaces, 
{\it e.g.}, $f_1 := f^{i_1}_{j_1}$
is a matrix which acts in the first vector space.
Additionally, we use  
incoming and outcoming indices, {\it e.g.},
$\Omega^{<1|} := \Omega^{i_1}$ and  $\gamma_{|1>} := \gamma_{j_1}$
denote a covector with one outcoming index and 
a vector with one incoming index respectively. 
Thus, in this notation, the matrix $f_1$ can be written as
$f_1 = f^{<1|}_{|1>}$.

\vspace{.5cm}

{\it Remark.}
Quantum Lie algebras defined by equations 
(\ref{int1})-(\ref{int1a}) generalize the
usual Lie (super-)algebras. Indeed in the non-deformed case, when
$$
\sigma^{mk}_{ij} = (-1)^{(m)(k)} \, \delta^m_j \, \delta^k_i
$$
is a super-permutation matrix (here $\sigma^2 =1$
and (\ref{int1a}) is fulfilled;
$(m)=0,1$ is the parity of a generator $\chi_m$), equations
(\ref{int1}) and (\ref{int2}) coincide with the defining 
relations and the Jacobi identities for Lie (super)-algebras.
Equation (\ref{int3}) is then equivalent to the 
$Z_2$-homogeneity condition $C^i_{jk} =0$ for $(i) \neq (j)+(k)$.
Equation (\ref{int3a}) follows from (\ref{int3}).

\vspace{0.5cm}

{\bf 2.} The exterior extension $\Gamma^{\wedge}$
of the quantum algebra $\Gamma$ (\ref{int1})
is obtained by adding new generators $\gamma_i$, $i=1, \dots, N$.
The generators $\gamma_i$ form a "generalized" wedge
algebra.
The definition of the wedge product of
 the elements $\gamma_i$ is
  \be
  \lb{13a}
 \gamma_{|1>} \wedge \gamma_{|2>} \dots \wedge \gamma_{|n>} =
 A_{1 \rightarrow n} \, \gamma_{|1>} \otimes \gamma_{|2>} \dots
 \otimes \gamma_{|n>} \; .
  \ee
Here the matrix operator $A_{1 \rightarrow n}$ is an analog of the
 antisymmetrizer of $n$-spaces. This operator can be defined
 inductively (see e.g. \cite{Gur})
\be
  \lb{anti} 
  A_{1\rightarrow n} = \left(
  {\bf 1} + \sum_{k=1}^{n-1} \, (-1)^{n-k} \,
 \sigma_{k \rightarrow n} \right) \,  A_{1\rightarrow n-1} 
\ee
 where, for $n >k$,
 $$
  \sigma_{k \rightarrow n} := \sigma_{k k+1} \, 
  \sigma_{k+1 k+2} \dots \sigma_{n-1 n}  \; .
  $$
Using the Yang-Baxter equation (\ref{int1a}) one can rewrite (\ref{anti})
in the following three equivalent forms 
$$
  A_{1\rightarrow n} = A_{1\rightarrow n-1} \, \left( {\bf 1}
 + \sum_{k=1}^{n-1} \, (-1)^{k} \, \sigma_{k+1 \leftarrow 1}
 \right)
 $$
 $$
 = \left( {\bf 1}
 + \sum_{k=1}^{n-1} \, (-1)^{k} \, \sigma_{k+1 \leftarrow 1}
 \right) \, A_{2\rightarrow n} =
 A_{2\rightarrow n}  \, \left({\bf 1} + \sum_{k=1}^{n-1} \, (-1)^{n-k} \,
 \sigma_{k \rightarrow n} \right) \; ,
 $$
 where
 $$
  \sigma_{n \leftarrow k} := \sigma_{n-1 n} \dots \sigma_{k+1 k+2} \, 
  \sigma_{k k+1} \; 
  $$
  for $n > k$.
  
 If the sequence of operators $A_{1 \rightarrow n}$ terminates
 at the step $n = h+1$ ($A_{1 \rightarrow h} \neq 0$ and
 $A_{1 \rightarrow n} = 0$ for $n >h$)
 then the number $h$ is called the height of 
 the operator $\sigma$.  

The cross-commutation relations between the generators
$\gamma_i$ and $\chi_j$ are:
\be
 \lb{app2.3} 
 \gamma_{|1>} \, \chi_{|2>} = 
 ( \sigma_{12} \,  \chi_{|1>}  +
 C^{<2|}_{|12>} ) \, \gamma_{|2>} \; . 
 \ee
The algebra $\Gamma^{\wedge}$ is graded by the degree in the
generators of $\gamma_i$.

\vspace{.5cm}

{\bf 3.} We further introduce a set of generators $\{ \Omega^i \}$, 
$i=1, \dots, N$,
canonically conjugated to the generators $\gamma_i$. The generators
$\Omega^i$ form a "wedge" algebra as well,
with the wedge product defined by 
  \be
  \lb{owedge}
   \Omega^{<r|} \wedge \Omega^{<r-1|} \wedge
   \dots \wedge \Omega^{<1|} = \Omega^{<r|} \otimes \Omega^{<r-1|}
   \otimes \dots \otimes \Omega^{<1|} \, A_{1 \rightarrow r} \; .
  \ee
Here operators $A_{1 \rightarrow n}$ are the same as in (\ref{anti}).
   
The commutation relations between $\Omega^i$ and $\gamma_j$ are
 \be
\lb{app2.1} 
\gamma_j \, \Omega^i = - \Omega^p \,
(\sigma^{-1})^{si}_{pj} \, \gamma_s + \delta^i_j  \, \Rightarrow
\, \gamma_{|2>} \, \Omega^{<2|} = - \Omega^{<1|} \,
\sigma^{-1}_{12} \, \gamma_{|1>} + I_2  \; .
 \ee
 
 Finally the commutation relations between $\Omega^i$ and $\chi_j$ are
 \be
 \lb{app2.2} 
 \chi_{|2>} \, \Omega^{<2|} =  \Omega^{<1|} \,
 ( \sigma_{12} \, \chi_{|1>} +  C^{<2|}_{|12>} ) \; . 
 \ee
We denote the algebra generated by $\{ \chi_i \}$,
$\{ \gamma_j \}$ and $\{ \Omega^k \}$ by $\Gamma^{\wedge} [ \Omega ]$.
The algebra $\Gamma^{\wedge}[\Omega]$ is graded by the rule:
$\mathrm{deg} \, (\gamma_i)=1$ and
$\mathrm{deg} \, (\Omega^i)=- 1$.
 
We shall need the following set of consequences
of the equation (\ref{app2.1}):
 $$
 \gamma_{|1>} \wedge \dots \wedge \gamma_{|r>} \, \Omega^{<r|} =
 (-1)^r \,
  \Omega^{<0|} \, \sigma^{-1}_{r \leftarrow 0} \,
  \gamma_{|0>} \wedge \dots \wedge \gamma_{|r-1>}
 $$
 \be
   \lb{app2.1a} 
 + (\sum_{k=1}^{r} \, (-1)^{r-k} \, 
 \sigma^{-1}_{r \leftarrow k} ) \, \gamma_{|1>} \wedge \dots \wedge
 \gamma_{|r-1>} \, ,
  \ee
 where $\sigma^{-1}_{r \leftarrow k} :=
  \sigma^{-1}_{k \, k+1} \dots \sigma^{-1}_{r-1 \, r}$ and
  $\sigma^{-1}_{r \leftarrow r} := {\bf 1}$ . 

\vspace{0.5cm}

{\bf 4.} The main result of the present paper is a recursive formula for the
BRST operator $Q$ which satisfies $Q^2 =0$. 

Such an operator endows
the algebra $\Gamma^{\wedge}$ with the
structure of the differential (chain) complex. To construct the differential
(starting with the operator $Q$)
one needs first to define the action of
the algebra $\Gamma^{\wedge}[ \Omega]$ on the algebra $\Gamma^\wedge$.
The elements $\chi_i$ and $\gamma_j$ act on $\Gamma^\wedge$
by the left multiplication.
To define the action of generators $\Omega^i$ on
$\Gamma^\wedge$ it suffices (due to relations (\ref{app2.1}) and
(\ref{app2.2})) to know $\Omega^i (1)$, where $1$ is the unit element
of the algebra $\Gamma^\wedge$. We set $\Omega^i (1)=0$.
The definition of the differential $d$ is given by
its action on an element
$\phi$ of the algebra $\Gamma^\wedge$, 
\be
\lb{dif}
d  \phi = [ Q , \phi ]_{_\pm} (1) \; ,
\ee
where $[,]_{_\pm}$ is the graded commutator. 

\vspace{0.5cm}
Now we are ready to formulate the main Proposition.

\vspace{0.5cm}
\noindent
{\bf Proposition.}
{\it The BRST operator $Q$ for the quantum algebra (\ref{int1})
has the following form
 \be
\lb{brst} Q = \Omega^i \, \chi_i + \sum^{h-1}_{r=1} \, Q_{(r)} \; ,
 \ee
where $h$ is the height of the operator $A_{1 \rightarrow n}$ (\ref{anti}).

Here the operators $Q_{(r)}$ are given by
 \be
\lb{qrr} Q_{(r)} =\Omega^{<r+1|} \, \Omega^{<r|} \dots
\Omega^{<1|} \, X^{<\tilde{1} \dots \tilde{r}|}_{|1 \dots r+1>} \,
\gamma_{|\tilde{1}>} \dots \gamma_{|\tilde{r}>}
 \ee
(the wedge product is implied); $X^{<1 \dots
r |}_{|1 \dots r+1>}$ are tensors which satisfy the following recurrent
relations
 \be
 \lb{app3.3}
A_{1 \rightarrow r+1} \,
 X^{< 1\dots r|}_{|1 \dots r+1>} \, A_{1 \rightarrow r } = 
A_{1 \rightarrow r+1} \, \left((-1)^{r} \, 
\sigma_{r+ 1 \leftarrow 1} - {\bf 1} \right) \, 
X^{< 2\dots r |}_{|2 \dots r+1>} \,
A_{2 \rightarrow
 r}  \; 
 \ee
with the initial condition $A_{12} \, X^{<0|}_{|12>} = - C^{<0|}_{|12>}$.
 }
 
\vspace{.5cm}
\noindent
{\bf Proof}. We have to verify the identity 
 \be
 \lb{Q2}
 Q^2 = ( \Omega^{<2|} \, \chi_{|2>})^2 +
 [ \Omega^{<2|} \, \chi_{|2>} \, , \, \sum_{r=1}^{h-1} \, Q_{(r)} ]_+ 
 + (\sum_{r=1}^{h-1} \, Q_{(r)})^2 = 0 \; .
 \ee 
Because of the lack of space we shall 
check a part of this identity which includes
the linear in $\chi$ terms only.

First of all we find (see (\ref{app2.2}))
$$
(\Omega^{<2|} \, \chi_{|2>})^2 =
\Omega^{<2|} \, \left( \Omega^{<1|} \,
 ( \sigma_{12} \, \chi_{|1>} +  C^{<2|}_{|12>} ) \, 
 \right) \,
\chi_{|2>} 
$$
\be
\lb{Q2a}
= \Omega^{<2|} \otimes \Omega^{<1|} \,
  \sigma_{12} \, C^{<2|}_{|12>} \, \chi_{|2>} +  
  \Omega^{<2|} \otimes \Omega^{<1|} \, ({\bf 1}-\sigma)_{12} \, C^{<2|}_{|12>} 
\, \chi_{|2>} 
\ee

$$
=
\Omega^{<2|} \otimes \Omega^{<1|} \, C^{<2|}_{|12>} \, \chi_{|2>} \; .
$$
Consider then the anticommutator
$[ \Omega^{<2|} \, \chi_{|2>} , Q_{(r)} ]_+ $ in which we commute
all $\chi_i$ to the right and extract
only the terms which are linear in the generators $\chi_i$:
 $$
[ Q_{(r)} , \, \Omega^{<r|} \chi_{|r>} ]_+ = \Omega^{<r+1|} \,
\chi_{|r+1>} Q_{(r)} + Q_{(r)} \, \Omega^{<r|} \, \chi_{|r>}
 $$
 
 $$ 
=\Omega^{<r+1|} \dots \Omega^{<0|} \, \left( \sigma_{r+1 \leftarrow 0}
+ (-1)^r {\bf 1} \right) \, X^{<1\dots r|}_{|1 \dots r+1>} \,
 \sigma^{-1}_{r \leftarrow 0} \, 
\gamma_{|0>} \dots \gamma_{|r-1>} \,
\chi_{|r>} 
 $$
 \be
 \lb{app2.4}
  + \Omega^{<r+1|} \dots
\Omega^{<1|} \, X^{< 1\dots r |}_{|1 \dots r+1>} \,
\left(\sum_{k=1}^{r} \, (-1)^{r-k} \,
\sigma^{-1}_{r \leftarrow k}  \right) \,
\gamma_{|1>} \dots \gamma_{|r-1>} \,
\chi_{|r>} + \dots \; 
 \ee
(dots denote the terms independent of $\chi_i$). Here  
eqs. (\ref{app2.2}), (\ref{app2.3}) and (\ref{app2.1a}) 
have been used. 

Equations (\ref{app2.4}) and (\ref{Q2a}) give the whole contribution
to the $\chi$-linear terms in $Q^2$
since $(\sum_{r=1}^{h-1} \, Q_{(r)})^2$ is independent of $\chi_i$.

The substitution of (\ref{Q2a}) and (\ref{app2.4})
produces the initial data $A_{12} X^{<0|}_{|12>} = -C^{<0|}_{|12>}$ and
recurrent relations
 \be
 \lb{app3.2}
 \begin{array}{c}
A_{1 \rightarrow r+1} \,
 X^{<1\dots r|}_{|1 \dots r+1>} \,
\left(\sum_{k=1}^{r} \, (-1)^{r-k} \,
\sigma^{-1}_{r \leftarrow k}  \right) \, A_{1 \rightarrow r-1} 
 \\ \\ 
= - A_{1 \rightarrow r+1} \, \left( \sigma_{r+1 \leftarrow 1} +
(-1)^{r-1} {\bf 1}\right) \, X^{<2\dots r|}_{|2 \dots r+1>}
\, \sigma^{-1}_{r \leftarrow 1} \, A_{1 \rightarrow r-1} \; ,
 \end{array}
 \ee
where the matrix operator $A_{1 \rightarrow r}$ is defined in
(\ref{anti}). These relations express coefficients
$X^{<1\dots r|}_{|1 \dots r+1>}$ via
$X^{<1\dots r-1|}_{|1 \dots r>}$.

Using an identity $\sigma^{-1}_{r \leftarrow 1} \, A_{1 \rightarrow
r-1}=A_{2\rightarrow r} \, \sigma^{-1}_{r \leftarrow 1}$ and
inductive relations (\ref{anti}) for the projectors $A_{1 \rightarrow
r}$ one can rewrite (\ref{app3.2}) in the form (\ref{app3.3}). $\bullet$ \\

\vspace{0.5cm}
 
 {\bf 5.} Comments.
 
 \vspace{.5cm}
 \noindent
 {\bf i.} For general $\sigma^{ij}_{kl}$ and $C^i_{jk}$
 it is rather difficult to solve equations (\ref{app3.3}) explicitly. 
 However for the case $\sigma^2 = {\bf 1}$ the main equations 
 (\ref{app3.3}) become simpler and the general solution for $Q$
 can be found. Indeed the
 relation (\ref{app3.3}) for $r=2$ gives
  $$
 A_{1 \rightarrow 3} \,
 X^{<12|}_{|123>} \, \left({\bf 1} - \sigma_{12}  \right) = A_{1
 \rightarrow 3} \, \left( \sigma_{23} \, \sigma_{12} -{\bf 1} \right) \,
 X^{<2|}_{|2 3>}  \; .
  $$
For $\sigma^2 =1$ 
we have $A_{1 \rightarrow 3} \, 
\left( \sigma_{23} \, \sigma_{12} -{\bf 1} \right) = 0$ and 
therefore $Q_{(r)} = 0$ for $r \geq 2$. Thus the BRST operator
(\ref{brst}) has the familiar form
$$
Q=\Omega^{<1|} \, \chi_{|1>} - \Omega^{<2|} \otimes \Omega^{<1|}  \,
C^{<1|}_{|1 2>} \, \gamma_{|1>} \; .
$$ 
In the case when the matrix $\sigma$ is the (super)-permutation matrix
the algebra $\Gamma^\wedge$ with the differential (\ref{dif})
becomes the standard complex for the Lie (super)-algebra $\Gamma$
(see {\it e.g.} \cite{Jac}).

In general, for $\sigma^2 \neq 1$, the sum in (\ref{brst})
will be limited only by the height $h$
of the operator $\sigma$.

Below we present an explicit form for $Q$ 
for the standard quantum deformation
$\Gamma = U_q(gl(N))$ of the universal enveloping algebra
of the Lie algebra $gl(N)$ ($\sigma^2 \neq 1$ in this case).

\vspace{0.5cm}
\noindent
{\bf ii.} When the algebra (\ref{int1}) is a Hopf
 algebra, the algebraic structure (\ref{13a}),
(\ref{app2.3}) -- (\ref{app2.2}) is
 related to the differential calculus on quantum groups
 (see \cite{Wor}, \cite{Schu}, \cite{RadVl}). 
 The BRST operator $Q$ given by (\ref{brst})
 generates the differential $d$ (introduced in \cite{Wor}) on the 
 algebra dual to $\Gamma^{\wedge}$.

\vspace{0.5cm}
 
{\bf 6.} Example. The BRST operator $Q$ for 
the quantum algebra $\Gamma =U_q(gl(N))$.

\vspace{0.5cm}

The quantum algebra $U_q(gl(N))$
is defined (as a Hopf algebra)
by the relations \cite{FRT}
\be
\label{92}
\R \, L^{\pm}_2 L^{\pm}_1  =  L^{\pm}_2 L^{\pm}_1 \, \R \,, \ \ \  
\R \, L^+_2 L^-_1  = L_2^- L^+_1 \, \R \,,  
\ee
\be
\label{93}
\Delta(L^{\pm})  =  L^{\pm}\otimes L^{\pm} , \ \ \
\varepsilon(L^{\pm}) = {\bf 1} , \ \ \
S(L^{\pm}) = (L^{\pm})^{-1} \; ,
\ee
where elements of the $N \times N$
matrices $(L^\pm)^i_j$ are generators of $U_q(gl(N))$;
the matrices $L^+$ and $L^-$ are 
respectively upper and lower triangular, 
their diagonal elements are related by $(L^+)^i_i \, (L^-)^i_i = 1$
for all $i$.
The matrix $\R$ is defined as
$\R := \R_{12}= P_{12} \, R_{12}$ ($P_{12}$ is the permutation matrix);
The matrix $R_{12}$ is the standard
Drinfeld-Jimbo $R$-matrix for $GL_q(N)$, 
$$
R_{12}=R^{i_{1},i_{2}}_{j_{1},j_{2}}=
\delta^{i_{1}}_{j_{1}} \delta^{i_{2}}_{j_{2}}(1+(q-1)\delta^{i_{1}i_{2}}) +
(q-q^{-1})\delta^{i_{1}}_{j_{2}} \delta^{i_{2}}_{j_{1}}
\Theta_{i_{1}i_{2}} \; , 
$$
where
$$
\Theta_{ij} = 
\left\{ 
\begin{array}{c}
1 \;\;\; {\rm if} \;\; i>j   \; ,  \\
0 \;\;\; {\rm if} \;\; i \leq j  \; .
\end{array}
\right.
$$
This $R$-matrix satisfies the Hecke condition 
$\R^2 = \lambda \, \R + {\bf 1}$, where $\lambda = (q-q^{-1})$
and $q$ is a parameter of deformation.

The generators of the algebra $\Gamma$ are defined by the formula  \cite{Jur},
\cite{Is}, \cite{RadVl}
\be
\label{102} 
\chi _{k}^{l}
=\frac{1}{\lambda}\,[\,(D^{-1})_{k}^{l}-(D^{-1})_{i}^{j}f_{kj}^{li}\,] \; .
\ee
Here 
$f_{kj}^{li}={L^-}_{k}^{i} S({L^+}_{j}^{l})$ and the
numerical matrix $D$ can be found by means of relations
$$
Tr_2 \R_{12} \Psi_{23} = P_{13} = Tr_2  \Psi_{12} \, \R_{23}  \; , \;\;\;
D_1 := Tr_2 \Psi_{12} \; \Rightarrow \; Tr_1(D^{-1}_1 \R^{-1}) = {\bf 1}_2 \; ,
$$
where $Tr_1$ and $Tr_2$ denote the traces over first and second spaces.

It is convenient to write down
the complete set of commutation relations for the exterior algebra
$\Gamma^{\wedge}[\Omega]$ in terms of generators 
$$
L^i_j = (L^+)^i_k S((L^-)^k_j) = \delta^i_j - \lambda \, 
S^{-1}(\chi^i_k) \, D^k_j \; , 
$$
$$
J^i_n = - S^{-1}(f^{ik}_{jl}) \, \gamma^l_k D^j_n \; , \;\;\;
\omega^i_j = \Omega^k_m \, f^{mi}_{kj} \; .
$$
The indices now are pairs of indices;
the roles of the elements $\chi_i$, $\gamma_j$ and $\Omega^k$
are played by the generators $\chi^i_j$, $\gamma^i_j$ and $\Omega^i_j$
respectively.

The commutation relations are \cite{SWZ}, \cite{Is}, \cite{RadVl}:
\be
\label{106} 
\omega_2 \R^{-1} \omega_2 \R = - \R^{-1} \omega_2 \R^{-1} \omega_2
\; , \;\;\;
\omega_2 \, \R \, L_2 \, \R = \R \, L_2 \, \R \, \omega_2 \; , 
\ee
\be
\label{109}
\omega_2 \, \R \, J_2 \R + \R \, J_2 \, \R \, \omega_2 = -\R \, ,
\;\;\;
L_2  \R \, L_2  \R = \R \, L_2 \R \, L_2 \; , 
\ee
\be
\label{113}
J_2 \R \, L_2 \R = \R \, L_2 \R \, J_2 \; , \;\;\;
J_2 \R \, J_2 \R = - \R^{-1} J_2 \R \, J_2 \, .
\ee

Now the construction of the BRST operator $Q$ is in order.
To begin we find the first term in the sum (\ref{brst}):
\be
\lb{114}
\Omega^k_m \, \chi^m_k =
  \frac{1}{\lambda} \, Tr_q \left( \omega \, (L- {\bf 1}) \right) 
 \; ,
\ee
where we have introduced the quantum trace $Tr_q (X) := Tr (D^{-1} X)$.
Then one can resolve the chain of the recurrent relations
(\ref{app3.3}) where we have to substitute the
expressions for the structure constants
$$
\sigma^{<^{jn}_{lq}|}_{|^{mi}_{p \, k}>}
=  R^{ju}_{sp} \, (R^{-1})^{sm}_{kr} \, (D^{-1})^f_o \,
R^{no}_{ut} \, D^t_l \, (R^{-1})^{ri}_{qf} \; ,
$$
$$
C^{<^q_p|}_{|^{im}_{jn}>}
= \delta^q_j \, \delta^i_n \, \delta^m_p - 
\sigma^{<^{qt}_{tp}|}_{|^{im}_{j \, n}>} \; ,
$$ 
and find the set of 
coefficients $X^{<1\dots r|}_{|1 \dots r+1>}$. After straightforward but
tiresome calculations one can obtain the following result:
$$
Q = Tr_q \left( \omega \, (L- {\bf 1})/ \lambda 
-    \omega \, L \, (\omega J) + \lambda \, \omega \, L \, (\omega J)^2 -
\lambda^2 \, \omega \, L \, (\omega J)^3 + \dots \right)
$$
$$
= Tr_q \left( \omega \, (L- {\bf 1})/ \lambda 
- \omega \, L \, (\omega J) \, ({\bf 1} + \lambda \omega J)^{-1} \right)
$$
\be
\lb{Qgl}
= - \frac{1}{\lambda} \, Tr_q ( \omega ) + \frac{1}{\lambda} \,
Tr_q \left( W  \right)
\; ,
\ee
where $W = \omega \, L\, ({\bf 1} + \lambda \, \omega J)^{-1}$
and the sum in the first line of (\ref{Qgl})
is limited by the requirement that monomials of $\omega$'s
of the order ${N^2+1}$ are equal to zero. 

One can check directly that the operator $Q$ 
given by (\ref{Qgl}) satisfies:
$$
 Q^2 = 0 \; ,
\;\;\; [ Q , \, L ] = 0 \; , \;\;\;
[ Q , \, J ]_+ = \frac{1}{\lambda} \, ({\bf 1}- L)  \; .
$$
To obtain these relations one has to use identities 
$$
Tr_q(X){\bf 1}_2 = Tr_{q1}(\R^{\pm 1} X_2 \R^{\mp 1})
$$ 
and relations
$$
\R \, W_2 \, \R^{-1} \, \omega_2 = - \omega_2 \, \R^{-1} \, W_2 \, \R  
\; , 
$$

$$
\R \, W_2 \, \R^{-1} \, W_2 = - W_2 \, \R^{-1} \, W_2 \, \R^{-1}  \; , 
\;\;\;
\R^{-1} \, W_2 \, \R \, L_2 = L_2 \, \R \, W_2 \, \R^{-1}  \; , 
$$

$$
J_2 \, \R \, W_2 \, \R^{-1} + \R^{-1} \, W_2 \, \R \, J_2 =
- L_2 \, ({\bf 1} + \lambda \, \omega J)^{-1}_2 \, \R^{-1} \, 
({\bf 1} + \lambda \, \omega J)_2 \; ,
$$
which follow from  (\ref{106})-(\ref{113}).

\vspace{.5cm}

{\it Remark.} The operator $Q$ given by (\ref{Qgl}) has the
correct classical limit for $q \rightarrow 1$ ($\lambda \rightarrow 0$,
$L \rightarrow {\bf 1} + \lambda \tilde{\chi}$, 
$\omega \rightarrow \tilde{\omega}$,
$J \rightarrow \tilde{\gamma}$)
$$
Q \rightarrow Q_{cl}=
Tr(\tilde{\omega} \tilde{\chi} + \tilde{\omega}^2 \, \tilde{\gamma} ) =
Tr(\tilde{\omega} \, X - \tilde{\omega} \, \tilde{\gamma} \, \tilde{\omega})
\; ,
$$
where $X := \tilde{\chi} + \tilde{\omega}
\, \tilde{\gamma} + \tilde{\gamma} \, \tilde{\omega}$ and
the classical algebra is
$$
[\tilde{\omega}_2, \, \tilde{\gamma_1} ]_+ = P_{12} \; , \;\;\;
[\tilde{\omega}_2, \, \tilde{\omega}_1 ]_+ = 0 
= [\tilde{\gamma}_2, \, \tilde{\gamma}_1 ]_+ \; , 
$$
$$
[ X_2 , \, X_1 ] = P_{12}(X_2 - X_1) \; , \;\;\;
[X_2, \, \tilde{\omega}_1 ] =  0  = [X_2, \, \tilde{\gamma}_1 ] \; .
$$

\subsection*{Acknowledgements.} 

We thank P. N. Pyatov, R. Stora and 
A. A. Vladimirov for valuable discussions. AI also thanks
K. Schmudgen and his team for useful comments and hospitality
at Leipzig University. This work was partially supported
by the RFBR grant 98-01-2033, the CNRS grant 
PICS-608 and the Votruba-Blokhintsev program. The work of AI was also supported
by the RFBR grant 00-01-00299.

\end{document}